\documentclass[12pt]{amsart}

\usepackage{geometry}
\usepackage{amsfonts, adjustbox, amssymb, amscd}
\numberwithin{equation}{section}

\usepackage{bm}
\usepackage{verbatim}
\usepackage{mathrsfs}
\usepackage{graphicx}
\usepackage{tikz-cd}
\usepackage{subcaption}
\usepackage{listings}
\usepackage{subfiles}
\usepackage[toc,page]{appendix}
\usepackage{mathtools}
\usepackage{comment}
\usepackage{enumerate}
\usepackage{enumitem}
\usepackage[all]{xy}
\usepackage{calligra}

\usepackage{graphicx}
\usepackage{appendix}
\usepackage{hyperref}
\hypersetup{
    colorlinks=true,
    citecolor=red,
    linkcolor=blue,
    filecolor=magenta,      
    urlcolor=red,
}
\newcommand{\Qq}{\mathbb{Q}}

\newcommand{\Rr}{\mathbb{R}}

\newcommand{\Exc}{\operatorname{Exc}}

\newcommand{\Mov}{{\operatorname{Mov}}}

\DeclareMathOperator{\HHom}{\mathscr{H}\text{\kern -3pt {\calligra\large om}}\,}

\newcommand{\Supp}{\operatorname{Supp}}

\newcommand{\mult}{\operatorname{mult}}


\newcounter{parentnumber}

\newtheorem{thm}{Theorem}[section]

\newtheorem{conj}[thm]{Conjecture}
\newtheorem{cor}[thm]{Corollary}
\newtheorem{lem}[thm]{Lemma}

\theoremstyle{definition}
\newtheorem{defn}[thm]{Definition}

\theoremstyle{definition}
\newtheorem{rem}[thm]{Remark}

\theoremstyle{definition}

\begin{document}

\title{Non-vanishing implies numerical dimension one abundance}
\dedicatory{To Professor Vyacheslav V. Shokurov on the occasion of his 75th birthday}

\author{Jihao Liu}
\author{Zheng Xu}

\date{\today}

\begin{abstract}
We show that the non-vanishing conjecture implies the abundance conjecture when $\nu\leq 1$. We also prove the abundance conjecture in dimension $\leq 5$ when $\kappa\geq 0$ and $\nu\leq 1$ unconditionally.
\end{abstract}

\address{Department of Mathematics, Peking University, No. 5 Yiheyuan Road, Haidian District, Beijing 100871, China}
\email{liujihao@math.pku.edu.cn}

\address{Beijing International Center for Mathematical Research,
Peking University, No. 5 Yiheyuan Road, Haidian District, Beijing 100871, China}
\email{zhengxu@pku.edu.cn}

\maketitle

\pagestyle{myheadings}\markboth{\hfill Jihao Liu and Zheng Xu \hfill}{\hfill Non-vanishing implies $\nu=1$ abundance\hfill}

\tableofcontents

\section{Introduction}\label{sec:Introduction}

We work over the field of complex numbers $\mathbb{C}$. 

The abundance conjecture predicts that Kodaira dimension $\kappa=\kappa(X)$ and the numerical dimension $\nu=\nu(X)$ (Definition \ref{defn: iitaka dimension}) of any smooth projective variety $X$ coincide:
\begin{conj}[Abundance]\label{conj: abundance}
Let $X$ be a smooth projective variety. Then $\kappa(X)=\nu(X)$.
\end{conj}
Conjecture \ref{conj: abundance} implies the existence of good minimal models or Mori fiber spaces for any smooth projective variety \cite[Theorem 4.3]{GL13} (see also \cite[Theorem 4.4]{Lai11}). 

The abundance conjecture is usually considered as (one of) the most important and difficult conjectures in birational geometry. In dimension $3$, the abundance conjecture was proven in \cite{Kaw92}, based on a series of earlier works \cite{Kaw85a,Kaw85b,Miy87,Miy88a,Miy88b} (see also \cite{KMM94,Fuj00}). The proof of the abundance conjecture for threefolds is based on a reduction via the Iitaka fibration to the case when either $\kappa=-\infty$ or $\kappa=0$ \cite{Kaw85a}, followed by a case-by-case treatment of $(\kappa,\nu)$: It is shown that $(\kappa,\nu)\neq(-\infty,0)$ \cite{Kaw85b}, $(\kappa,\nu)\neq(-\infty,2)$ \cite{Miy87}, $(\kappa,\nu)\neq(-\infty,1)$ \cite{Miy88a}, $(\kappa,\nu)\neq(0,1)$ \cite{Miy88b} (see \cite{Kaw92} for an alternative proof), and finally $(\kappa,\nu)\neq(0,2)$ \cite{Kaw92}, which concluded the proof of the abundance conjecture for threefolds.

In dimensions $\geq 4$, thanks to \cite{GL13,Lai11} and \cite{BCHM10}, we can still reduce the abundance conjecture to the case when $\kappa=-\infty$ or $\kappa=0$. Moreover, one can still show that $(\kappa,\nu)\neq(-\infty,0)$ \cite{Nak04} (see also \cite{CP11,Dru11,CKP12,Kaw13}) in any dimension, which implies the abundance conjecture when $\nu=0$. However, even when $\dim X=4$, no other non-trivial case of $(\kappa,\nu)$ could be excluded despite many partial results (cf. \cite{DHP13,GM17,LP18,LP20b,IMM24,Laz24}).

\medskip

\noindent\textbf{Main theorems.} The goal of this paper is to systematically study the case when $(\kappa,\nu)=(0,1)$ in higher dimensions (dimension $\geq 4$). As the first result of our paper, we resolve the abundance conjecture when $(\kappa,\nu)=(0,1)$ and $\dim X\leq 5$:

\begin{thm}\label{thm: main dim 5}
    Let $X$ be a smooth projective variety of dimension $\leq 5$ such that $\kappa(X)\geq 0$ and $\nu(X)\leq 1$. Then $X$ has a good minimal model. In particular, $\kappa(X)=\nu(X)$.
\end{thm}

Combining Theorem \ref{thm: main dim 5} with \cite[Theorem 6.7]{LP18}, we obtain the following:

\begin{cor}\label{cor: lp18}
    Let $X$ be a smooth projective variety of dimension $\leq 4$ such that $\chi(\mathcal{O}_X)\not=0$ and $\nu(X)\leq 1$.  Then $X$ has a good minimal model. In particular, $\kappa(X)=\nu(X)$.
\end{cor}

The second main result of our paper shows that the non-vanishing conjecture for smooth projective varieties (Conjecture \ref{conj: non-vanishing smooth}) implies the abundance conjecture when $\nu\leq 1$:

\begin{thm}\label{thm: main nv}
    Assume the non-vanishing conjecture for smooth projective varieties in dimension $d$. Then any smooth projective variety $X$ of dimension $\leq d$ with $\nu(X)\leq 1$ has a good minimal model or a Mori fiber space. In particular, $\kappa(X)=\nu(X)$.
\end{thm}

\begin{conj}[Non-vanishing for smooth projective varieties]\label{conj: non-vanishing smooth}
    Let $X$ be a smooth projective variety. If $\nu(X)\geq 0$, then $\kappa(X)\geq 0$.
\end{conj}

We remark that both Theorems \ref{thm: main dim 5} and \ref{thm: main nv} also hold if we consider log canonical (lc) pairs $(X,B)$ instead of smooth varieties $X$. See Theorems \ref{thm: lc abundance dim 5 nu 1} and \ref{thm: abundance lc case assuming non-vanishing} for details.

The reader may want to compare Theorem \ref{thm: main nv} with the main result of \cite{DHP13}, which indicates that the abundance conjecture follows from the non-vanishing conjecture and the dlt extension conjecture \cite[Conjecture 1.3]{DHP13}. In Theorem \ref{thm: main nv}, we remove the assumption of the dlt extension conjecture with the price of assuming that $\nu(X)\leq 1$. However, unlike \cite{DHP13} or any other recent progress toward the abundance conjecture that heavily relies on analytic methods, our approach is much more algebraic. The main ingredients in our proof are:
\begin{itemize}
\item[(i)] Kawamata's alternative approach to the threefold case \cite[Section 4]{Kaw92}.
\item[(ii)] The Du Bois property of slc singularities \cite{KK10,Kol13}.
\item[(iii)] Recent progress in the lc minimal model program (cf. \cite{HH20,LT22,MZ23}).
\end{itemize}

We note that (i) and (ii) were also considered in Kawamata's withdrawn note \cite{Kaw11} (see Remark \ref{rem: kaw11's idea} for a detailed explanation), while (iii) is a completely new input. In fact, thanks to recent advances in the minimal model program, both Theorems \ref{thm: main dim 5} and \ref{thm: main nv} can be reduced to the following more technical but stronger result, which shows that a special case of the special termination conjecture (Conjecture \ref{conj: stof with scaling weak}) implies the abundance conjecture when $(\kappa,\nu)=(0,1)$:

\begin{thm}\label{thm: main sp tof}
Assume Conjecture \ref{conj: stof with scaling weak} in dimension $\leq d$. Let $X$ be a smooth projective variety of dimension $\leq d$ such that $\kappa(X)\geq 0$ and $\nu(X)\leq 1$. Then $X$ has a good minimal model. In particular, $\kappa(X)=\nu(X)$.
\end{thm}

\begin{conj}[Special termination of some MMP for projective effective pairs]\label{conj: stof with scaling weak}
    Let $(X,B)$ be a projective $\mathbb{Q}$-factorial effective dlt pair. Then there exists a $(K_X+B)$-MMP that terminates near the divisorial part of the image of $\lfloor B\rfloor$.
\end{conj}

It is worth mentioning that Conjecture \ref{conj: stof with scaling weak} is implied by most standard conjectures in the minimal model program in \emph{\textbf{lower}} dimensions, such as the existence of minimal models (cf. \cite{LT22}), termination of MMP with scaling (cf. \cite{Bir10}), and non-vanishing (see Theorem \ref{thm: 1.5 low imply 1.7}).

\medskip

We remark that Theorem \ref{thm: main sp tof} also holds if we consider log canonical pairs $(X,B)$ instead of smooth varieties $X$. See Theorem \ref{thm: projective case} for details. All our main theorems still hold if we consider klt pairs $(X,B)$ associated with projective morphisms $X\rightarrow U$, and the condition ``$(\kappa,\nu)=(0,1)$" can be weakened to ``$\nu-\kappa\leq 1$". See Theorem \ref{thm: relative abundance nu=1} for a precise statement. For lc pairs, however, due to technical difficulties, we cannot generalize our statements to the relative case; see Remark \ref{rem: why lc not ok} for an explanation.

\medskip

\noindent\textbf{Acknowledgements.} The authors would like to thank Christopher D. Hacon, Jingjun Han, Junpeng Jiao, Vladimir Lazić, Fanjun Meng, Wenhao Ou, Nikolaos Tsakanikas, Weimufei Wu, and Lingyao Xie for useful discussions. The first author would like to thank Vyacheslav V. Shokurov for his insight; he mentioned to the author and/or Jingjun Han several years ago that Theorem \ref{thm: main nv} is approachable. Although the first author misunderstood and thought it referred to the $\nu=\dim X-1$ case, Shokurov's insight remains invaluable to him. The authors would like to thank Frédéric Campana, whose visit to Peking University and a lunch inspired the collaboration between the authors. The first author is supported by the National Key R\&D Program of China (\#2024YFA1014400). A part of the work was done when the second author visited Jakub Witaszek at Princeton University, and the second author would like to thank him for his hospitality. The second author would like to thank Zhiyu Tian for his constant support.

\section{Preliminaries}

We will adopt the standard notation and definitions on MMP in \cite{KM98,BCHM10,Sho92} and use them freely. 

\subsection{Pairs}

\begin{defn}
A \emph{pair} $(X,B)/U$ consists of a projective morphism $X\rightarrow U$ between normal quasi-projective varieties and an  $\Rr$-divisor $B\geq 0$ on $X$ such that $K_X+B$ is $\mathbb R$-Cartier. If $U=\{pt\}$, we usually drop $U$ and say that $(X,B)$ is \emph{projective}. If $U$ is not important, we may also drop $U$. For any prime divisor $E$ and $\mathbb R$-divisor $D$ on $X$, we define $\mult_{E}D$ to be the \emph{multiplicity} of $E$ along $D$.  Let $h: W\to X$
	be any log resolution of $(X,\Supp B)$ and let
	$$K_W+B_W:=h^*(K_X+B).$$
	The \emph{log discrepancy} of a prime divisor $D$ on $W$ with respect to $(X,B)$ is $1-\mult_{D}B_W$ and it is denoted by $a(D,X,B).$
	
We say that $(X,B)$ is \emph{lc} (resp. \emph{klt}) if $a(D,X,B)\geq 0$ (resp. $>0$) for every log resolution $h: W\to X$ as above and every prime divisor $D$ on $W$. An \emph{lc place} of $(X,B)$ is a prime divisor $E$ over $X$ such that $a(E,X,B)=0$. We say that $(X,B)$ is \emph{dlt} if there exists a log resolution of $(X,B)$ which only extracts divisors that are not lc places of $(X,B)$.

A projective pair $(X,B)$ is called \emph{effective} if $K_X+B\sim_{\mathbb R}D\geq 0$ for some $D$.
\end{defn}

\begin{defn}
    Two pairs $(X,B)$ and $(X',B')$ are called \emph{crepant} to each other if there exist birational morphisms $p: W\rightarrow X$ and $q: W\rightarrow X'$ such that $p^*(K_X+B)=q^*(K_{X'}+B')$.
\end{defn}

\subsection{Iitaka dimensions}
\begin{defn}\label{defn: iitaka dimension}
Let $X$ be a normal projective variety and $D$ an $\Rr$-divisor on $X$. The \emph{Iitaka dimension} $\kappa(D)$ (resp. \emph{numerical Iitaka dimension} $\kappa_{\sigma}(D)$) of $D$ is defined in the following way. If $|\lfloor mD\rfloor|\not=\emptyset$ for some positive integer $m$ (resp. $D$ is pseudo-effective), then
$$\kappa(D):=\max\left\{k\in\mathbb N\middle| \underset{m\rightarrow+\infty}{\lim\sup}\frac{\dim H^0(X,\lfloor mD\rfloor)}{m^k}>0\right\}$$
$$\left(\text{resp. }\kappa_{\sigma}(D):=\max\left\{k\in\mathbb N\middle| A\text{ is Cartier}, \underset{m\rightarrow+\infty}{\lim\sup}\frac{\dim H^0(X,\lfloor mD\rfloor+A)}{m^k}>0\right\}\right).$$
Otherwise, we let $\kappa(D):=-\infty$ (resp. $\kappa_{\sigma}(D):=-\infty$). If $|D|_{\mathbb R}\not=\emptyset$, then we define $\kappa_{\iota}(D):=\kappa(D')$ for some $D'\in |D|_{\mathbb R}$. Otherwise, we define $\kappa_{\iota}(D):=-\infty$. $\kappa_{\iota}$ is well-defined by \cite[Section 2]{Cho08}. We define $\kappa(X):=\kappa(K_X)$ as the \emph{Kodaira dimension} of $X$, and define $\nu(X):=\kappa_{\sigma}(K_X)$ as the \emph{numerical dimension} of $X$.

Let $\pi: X\rightarrow U$ be a projective morphism from between normal quasi-projective varieties and $D$ an $\Rr$-divisor on $X$. Let $X\xrightarrow{f} T\rightarrow\pi(U)$ be the Stein factorization and let $F$ be a general fiber of $f$. We define $\kappa(X/U,D):=\kappa(D|_F),\kappa_{\iota}(X/U,D):=\kappa_{\iota}(D|_F)$, and $\kappa_{\sigma}(X/U,D):=\kappa_{\sigma}(D|_F)$ as the \emph{Iitaka dimension}, \emph{invariant Iitaka dimension}, and \emph{numerical Iitaka dimension} of $D$ over $U$ respectively. Here by convention, if $\dim F=0$, then we define $\kappa(D|_F):=\kappa_{\iota}(D|_F):=\kappa_{\sigma}(D|_F):=0$.

For any non-normal projective variety $X$ with normalization $\nu: \tilde X\rightarrow X$ and $\mathbb R$-Cartier $\mathbb R$-divisor $D$ on $X$, we define $\kappa(D):=\kappa(\nu^*D)$, $\kappa_{\iota}(D):=\kappa_{\iota}(\nu^*D)$, and $\kappa_{\sigma}(D):=\kappa_{\sigma}(\nu^*D)$.

We refer the reader to \cite[Chapters II,V]{Nak04}, \cite[Section 2]{Cho08}, and \cite[Section 2]{HH20} for basic properties of Iitaka dimensions.
\end{defn}

\begin{lem}[{cf. \cite[Lemma 2.3]{LX23}}]\label{lem: property of numerical and Iitaka dimension} Let $\pi: X\rightarrow U$ be a projective morphism between normal quasi-projective varieties and $D$ an $\Rr$-Cartier $\Rr$-divisor on $X$. Then:
\begin{enumerate}
    \item Let $f: Y\rightarrow X$ be a surjective birational morphism and $D_Y$ an $\Rr$-Cartier $\Rr$-divisor on $Y$ such that $D_Y=f^*D+E$ for some $f$-exceptional $\Rr$-divisor $E\geq 0$. Then $\kappa_{\sigma}(Y/U,D_Y)=\kappa_{\sigma}(X/U,D)$ and $\kappa_{\iota}(Y/U,D_Y)=\kappa_{\iota}(X/U,D)$. 
    \item Let $g: Z\rightarrow X$ be a projective morphism from a normal quasi-projective variety $Z$. Then $\kappa_{\sigma}(Z/U,g^*D)=\kappa_{\sigma}(X/U,D)$ and $\kappa_{\iota}(Z/U,g^*D)=\kappa_{\iota}(X/U,D)$.   
    \item  Let $\phi: X\dashrightarrow X'$ be a partial $D$-MMP$/U$ and let $D':=\phi_*D$. Then $\kappa_{\sigma}(X/U,D)=\kappa_{\sigma}(X'/U,D')$ and $\kappa_{\iota}(X/U,D)=\kappa_{\iota}(X'/U,D')$
\end{enumerate}
\end{lem}

\begin{lem}\label{lem: compare kappa}
Let $X$ be a normal projective variety. Let Let $A\geq B\geq C$ be three $\mathbb R$-divisors on $X$ such that $C$ is pseudo-effective (resp. effective) and $\Supp(B-C)=\Supp(A-C)$. Then $\kappa_{\sigma}(B)=\kappa_{\sigma}(A)$ (resp. $\kappa_{\iota}(B)=\kappa_{\iota}(A)$).
\end{lem}
\begin{proof}
    It is obvious that $\kappa_{\sigma}(B)\leq \kappa_{\sigma}(A)$ and $\kappa_{\iota}(B)\leq \kappa_{\iota}(A)$. Since $\Supp(B-C)=\Supp(A-C)$, we have $C+\lambda (B-C)\geq A$ for some $\lambda\geq 1$, hence $\lambda B\geq A+(\lambda -1)C$, and so $\kappa_{\sigma}(B)\geq \kappa_{\sigma}(A)$ (resp. $\kappa_{\iota}(B)\geq \kappa_{\iota}(A)$). The lemma follows.
\end{proof}

\subsection{Models}

\begin{defn}
    Let $(X,B)$ be an lc pair. A \emph{$\mathbb Q$-factorial dlt modification} of $(X,B)$ is a birational morphism $h: X'\rightarrow X$ such that $(X',B':=h^{-1}_*B+\Exc(h))$ is $\mathbb Q$-factorial dlt and $K_{X'}+B'=h^*(K_X+B)$. $(X',B')$ is called a \emph{$\mathbb Q$-factorial dlt model} of $(X,B)$. Existence of $\mathbb Q$-factorial dlt modification is well-known, cf. \cite[Theorem 3.1]{KK10}.
\end{defn}

\begin{defn}[Log birational model]\label{defn: log birational model}
  Let $(X,B)/U$ be an lc pair, $\phi: X\dashrightarrow X'$ a birational map over $U$ and $E:=\Exc(\phi^{-1})$ the reduced $\phi^{-1}$-exceptional divisor. We let $B':=\phi_*B+E$
and say that $(X',B')/U$ is a \emph{log birational model} of $(X,B)/U$.
\end{defn}

\begin{defn}[Minimal models]\label{defn: minimal model}
    Let $(X,B)/U$ be an lc pair and $(X',B')/U$ a log birational model of $(X,B)/U$ such that $K_{X'}+B'$ is nef$/U$.
    \begin{enumerate}
        \item We say that $(X',B')/U$ is a \emph{minimal model in the sense of Birkar-Shokurov} of $(X,B)/U$, if for any prime divisor $D$ on $X$ which is exceptional over $X'$, $$a(D,X,B)<a(D,X',B').$$
        \item We say that $(X',B')/U$ is \emph{good minimal model in the sense of Birkar-Shokurov} of $(X,B)/U$ if it is a bs-minimal model of $(X,B)/U$ and $K_{X'}+B'$ is semi-ample$/U$. 
        \item We say that $(X',B')/U$ is a \emph{log minimal model} of $(X,B)/U$ if it is a minimal model in the sense of Birkar-Shokurov of $(X,B)/U$ and $(X',B')$ is $\Qq$-factorial dlt.
        \item We say that $(X',B')/U$ is a \emph{good log minimal model} of $(X,B)/U$ if it is a log minimal model of $(X,B)/U$ and $K_{X'}+B'$ is semi-ample$/U$. 
\end{enumerate}
If the induced birational map $X\dashrightarrow X'$ does not extract any divisor, then we remove the phrase ``in the sense of Birkar-Shokurov" in the descriptions of (1-2).
\end{defn}

It is important to remark that the definition of minimal models in references sometimes do not coincide. Fortunately for us, we have the following lemmas which allows us to avoid any potential issues.

\begin{lem}\label{lem: equivalence of log minimal model}
    Let $(X,B)/U$ be an lc pair. The followings are equivalent.
    \begin{enumerate}
        \item $(X,B)/U$ has a log minimal model (resp. good log minimal model)
        \item $(X,B)/U$ has a minimal model (resp. good minimal model)
        \item We may run a $(K_X+B)$-MMP$/U$ with scaling of an ample divisor which terminates with a minimal model (resp. good minimal model) of $(X,B)/U$.
    \end{enumerate}
\end{lem}
\begin{proof}
   The minimal model part part follows from \cite[Theorem 1.7]{HH20}. Combining with \cite[Remark 2.7]{Bir12}, we get the good minimal model part.
\end{proof}

The following several results on the minimal model program, some of which have been established in recent years, are important for the proof of our main theorems.

\begin{lem}\label{lem: equivalence existence of model}
    Let $(X,B)/U$ and $(X',B')/U$ be lc pairs with birational morphism $h: X'\rightarrow X$ such that $K_{X'}+B'=h^*(K_X+B)+E$ for some $E\geq 0$ that is exceptional$/X$. Then $(X,B)/U$ has a minimal model (resp. good minimal model) if and only if $(X',B')/U$ has a minimal model (resp. good minimal model).
\end{lem}
\begin{proof}
It follows from \cite[Lemma 2.15]{Has19a} and Lemma \ref{lem: equivalence of log minimal model}.
\end{proof}

\subsection{Preliminaries on MMP}

\begin{thm}[{\cite{Dru11,Gon11,Nak04}}]\label{thm: num 0 abundance}
    Let $(X,B)$ be a projective lc pair such that $\kappa_{\sigma}(K_X+B)=0$. Then:
    \begin{enumerate}
        \item $\kappa_{\iota}(K_X+B)=0$.
        \item If $K_X+B\equiv 0$, then $K_X+B\sim_{\mathbb R}0$.
        \item $(X,B)$ has a good minimal model.
    \end{enumerate}
\end{thm}
\begin{proof}
    (1) follows from \cite[Theorem 6.1]{Gon11} and (2) follows from \cite[Lemma 6.1]{Gon11}. Let $(X',B')$ be a $\mathbb Q$-factorial dlt model of $(X,B)$. By \cite[Theorem 1.1]{Gon11}, $(X',B')$ has a minimal model. By Lemma \ref{lem: equivalence existence of model}, $(X,B)$ has a minimal model $(X_{\min},B_{\min})$. By (2), $(X_{\min},B_{\min})$ is a good minimal model of $(X,B)$, which implies (3).
\end{proof}

\begin{lem}[{\cite[Lemma 2.13]{HH20}, \cite[Theorem 4.3]{GL13}}]\label{lem: hh20 2.13}
    Let $(X,B)/U$ be a klt pair such that $K_X+B$ is pseudo-effective$/U$ and abundant$/U$. Then $(X,B)/U$ has a good minimal model.
\end{lem}

\begin{thm}[{\cite[Theorem 1.4]{MZ23}}]\label{thm: mz23 1.4}
  Let $(X,B)/U$ and $(X,B')/U$ be two lc pairs such that $B\geq B'\geq 0$ and $\Supp B=\Supp B'$. If $(X,B)/U$ has a good minimal model and $K_X+B'$ is pseudo-effective$/U$, then $(X,B')/U$ has a good minimal model.
\end{thm}

\begin{lem}[{cf. \cite[Lemma 3.21]{HL22}}]\label{lem: hl22 3.21}
Let $(X,B)/U$ be a $\mathbb Q$-factorial lc pair such that $X$ is klt and $K_X+B$ is nef$/U$. Let $D$ be an $\mathbb R$-divisor on $X$ such that $(X,B+D)$ is an lc pair. Then for any $0<\epsilon\ll 1$, any sequence of steps of a $(K_X+B+\epsilon D)$-MMP$/U$ is $(K_X+B)$-trivial. 
\end{lem}

\begin{thm}\label{thm: 1.5 low imply 1.7}
 Conjecture \ref{conj: non-vanishing smooth} in dimension $d-1$ implies Conjecture \ref{conj: stof with scaling weak} in dimension $\leq d$.
\end{thm}
\begin{proof}
    By \cite[Theorem 1.4]{Has18},  Conjecture \ref{conj: non-vanishing smooth} in dimension $d-1$ implies the existence of minimal models for any smooth projective variety $X$ of dimension $\leq d-1$ such that $K_X$ is pseudo-effective. By \cite[Theorem B]{LT22}, this implies the existence of minimal models for any projective lc pair $(X,B)$ of dimension $\leq d$ such that $\kappa_{\iota}(K_X+B)\geq 0$. By \cite[Theorem 4.1]{Bir12}, for any projective $\mathbb Q$-factorial dlt pair $(X,B)$ of dimension $\leq d$, we may run a $(K_X+B)$-MMP with scaling of an ample divisor which terminates. In particular, it terminates near the divisorial part of the image of $\lfloor B\rfloor$.
\end{proof}

\subsection{Real coefficients}

We prove the $\mathbb R$-divisor version of some well-known results on $\mathbb Q$-divisors.

\begin{defn}
Let $m$ be a positive integer and $\bm{v}\in\mathbb R^m$. The \emph{rational envelope} of $\bm{v}$ is the minimal rational affine subspace of $\mathbb R^m$ which contains $\bm{v}$. For example, if $m=2$ and $\bm{v}=\left(\frac{\sqrt{2}}{2},1-\frac{\sqrt{2}}{2}\right)$, then the rational envelope of $\bm{v}$ is $(x_1+x_2=1)\subset\mathbb R^2_{x_1x_2}$.
\end{defn}

\begin{lem}\label{lem: ambro05 r coefficient}
    Let $(X,B)$ be a klt pair and let $f: X\rightarrow Z$ be a contraction such that $K_X+B\sim_{\mathbb R,Z}0$. Then there exists a klt pair $(Z,B_Z)$ such that $K_X+B\sim_{\mathbb R}f^*(K_Z+B_Z)$.
\end{lem}
\begin{proof}    
By \cite[Lemma 5.3, Theorem 5.6]{HLS24} there exist positive real numbers $a_1,\dots,a_k$ and $\mathbb Q$-divisors $B_1,\dots,B_k$, such that $\sum_{i=1}^ka_i=1$, $\sum_{i=1}^ka_iB_i=B$, each $(X,B_i)$ is klt, and $K_{X}+B_i\sim_{\mathbb Q,Z}0$. By \cite[Theorem 0.2]{Amb05} there exist klt $\mathbb Q$-pairs $(Z,B_{Z,i})$ such that $K_{X}+B_i\sim_{\mathbb Q}f^*(K_Z+B_{Z,i})$. We may let $B_Z:=\sum_{i=1}^ka_iB_{Z,i}$.
\end{proof}

\begin{lem}\label{lem: ckp12 real}
    Let $(X,B)$ be a projective lc pair such that $K_X+B\equiv D\geq 0$ for some $\mathbb R$-divisor $D$. Then $\kappa_{\iota}(K_X+B)\geq\kappa(D)\geq 0$.
\end{lem}
\begin{proof}
    We write $B=\sum_{i=1}^m b_iB_i$ and $D=\sum_{j=1}^n c_jD_j$ where $B_i$ are the irreducible components of $B$ and $D_j$ are the irreducible components of $D$. Let $V$ be the rational envelope of $\bm{v}_0:=(b_1,\dots,b_m,c_1,\dots,c_n)$ in $\mathbb R^{m+n}$. By \cite[Lemma 5.3,Theorem 5.6]{HLS24}, for any $\bm{v}=(b_1',\dots,b_m',c_1',\dots,c_n')\in V$ such that $||\bm{v}-\bm{v}_0||_{\infty}\ll 1$, we have that $(X,B(\bm{v}):=\sum_{i=1}^mb_i'B_i)$ is lc, $D(\bm{v}):=\sum_{i=1}^nc_i'D_i\geq 0$, and $K_X+B(\bm{v})\equiv D(\bm{v})$. Pick rational vectors $\bm{v}_1,\dots,\bm{v}_k\in V$ such that $||\bm{v}_i-\bm{v}_0||_{\infty}\ll 1$ for any $i$ such that $\bm{v}_0$ is contained in the interior of the convex hull spanned by $\bm{v}_1,\dots,\bm{v}_k$ and let $a_1,\dots,a_k$ be positive real numbers such that $\sum_{i=1}^ka_i=1$ and $\sum_{i=1}^ka_i\bm{v}_i=\bm{v}_0$. By \cite[Corollary 3.2]{CKP12}, $K_X+B(\bm{v}_i)\sim_{\mathbb R}L_i\geq 0$ for some $L_i$ such that $\kappa(L_i)\geq\kappa(D(\bm{v}_i)$. We have
    $$K_X+B\sim_{\mathbb R}\sum_{i=1}^ka_iL_i\geq 0$$
    and
    $$\kappa_{\iota}(K_X+B)=\kappa\left(\sum_{i=1}^ka_iL_i\right)\geq\kappa\left(\sum_{i=1}^ka_iD(\bm{v}_i)\right)=\kappa(D).$$
    The lemma follows.
\end{proof}

\begin{conj}[Log abundance]\label{conj: log abundance}
    Let $(X,B)$ be a projective lc pair such that $K_X+B$ is pseudo-effective. Then $(X,B)$ has a good minimal model.
\end{conj}

\begin{thm}\label{thm: abundance dim 3}
  Conjecture \ref{conj: log abundance} holds in dimension $\leq 3$.  
\end{thm}
\begin{proof}
This should be well-known as an easy consequence of \cite[1.1 Theorem]{KMM94} and perturbation. Let $(X,B)$ be a projective lc pair of dimension $3$ such that $K_X+B$ is pseudo-effective. Possibly replacing $(X,B)$ with a $\mathbb Q$-factorial dlt model, we may assume that $(X,B)$ is $\mathbb Q$-factorial dlt. By \cite[2.3 Theorem]{Sho96} we may assume that $K_X+B$ is nef. By \cite[Remark 3.1, Proposition 3.2]{Bir11}, we may write $K_X+B=\sum_{i=1}^k a_i(K_X+B_i)$ where each $a_i>0$, $\sum_{i=1}^ka_i=1$, each $(X,B_i)$ is an lc $\mathbb Q$-pair, and $K_X+B_i$ is nef for each $i$. By \cite[1.1 Theorem]{KMM94}, $K_X+B_i$ is semi-ample for each $i$, so $K_X+B$ is semi-ample.
\end{proof}

\subsection{Numerical dimension one}

Nef divisors with numerical dimension one have some special properties compared to those with larger numerical dimensions. We need the following result:

\begin{lem}\label{lem: restrictio of nu=1 is nu=0}
    Let $D=\sum a_iD_i\geq 0$ be a nef $\mathbb R$-divisor on a normal projective variety $X$ such that $\kappa_{\sigma}(D)=1$, where each $a_i>0$ and $D_i$ are prime divisors. Then $D|_{D_i}\equiv 0$ for any $i$.
\end{lem}
\begin{proof}
    Let $d:=\dim X$ and let $H$ be an ample divisor on $X$. Then 
    $$D\cdot D_i\cdot H^{d-2}\geq 0$$
    for each $i$, and
    $$\sum_{i=1}^na_i(D\cdot D_i\cdot H^{d-2})=D^2\cdot H^{d-2}=0.$$
    Thus $D\cdot D_i\cdot H^{d-2}=0$ for each $i$, hence
    $$D|_{D_i}\cdot (H|_{D_i})^{d-2}=0,$$
    so $\kappa_{\sigma}(D|_{D_i})=0$, and so $D|_{D_i}\equiv 0$. 
\end{proof}

The following result is an $\mathbb R$-divisor version of \cite[Theorem 6.1]{LP18} and is slightly strengthened.
\begin{thm}\label{thm: lp18 r coefficient}
    Let $X$ be a projective variety and $L$ a nef $\mathbb R$-divisor on $X$ such that $\kappa_{\sigma}(L)=1$. Let $D\geq 0$ be an $\mathbb R$-Cartier $\mathbb R$-divisor on $X$ and $F$ a pseudo-effective $\mathbb R$-Cartier $\mathbb R$-divisor on $X$, such that $D+F\equiv L$. Pick $s\in (0,1]$. Then there exists an $\mathbb R$-divisor $E\geq 0$ on $X$ such that $D+sF\equiv E$ and $\kappa(E)\geq\kappa(D)$.
\end{thm}
\begin{proof}
    Let $h: X'\rightarrow X$ be a resolution of $X$. Possibly replacing $D,F,L$ with their pullbacks to $X'$, we may assume that $X$ is smooth. Let $P$ and $N$ be the positive and negative part of the Nakayama-Zariski decomposition of $F$ respectively, then we have $F=P+N$. If $P\equiv 0$, then $E:=D+sN$ satisfies our requirements, so we may assume that $P\not\equiv 0$. Let $S$ be a surface in $Y$ cut out by $\dim X-2$ general hyperplane sections. Since $P\in\overline{\Mov}(X)$, $P|_S\in\overline{\Mov}(S)$. Since $S$ is a surface, $P|_S$ is nef. Thus $(P|_S)^2\geq 0$. Since $\kappa_{\sigma}(L)=1$, we have
$$0=(L|_S)^2=L|_S\cdot P_S+L|_S\cdot N|_S+L|_S\cdot D|_S.$$
Since $L|_S$ is nef, 
$$L|_S\cdot P_S=L|_S\cdot N|_S=L|_S\cdot D|_S=0.$$
The Hodge index theorem implies $(P|_S)^2\leq 0$, hence $(P|_S)^2=0$. By \cite[Lemma 2.2]{LP18}, $P|_S\equiv\lambda L|_S$ for some real number $\lambda>0$, and hence $P\equiv\lambda L$ by the Lefschetz hyperplane section theorem. Since $D\not=0$, $\lambda<1$. We have
$$D+\lambda L+N\equiv D+P+N=D+F\equiv L,$$
so we have
$$\frac{1}{1-\lambda}(D+N)\equiv L\equiv D+F$$
and so
$$D+sF=s(D+F)+(1-s)D\equiv\frac{s}{1-\lambda}(D+N)+(1-s)D\geq\frac{s}{1-\lambda}D.$$
We may take 
$$E:=\frac{s}{1-\lambda}(D+N)+(1-s)D.$$
\end{proof}

\section{A special case}

In this section we shall consider semi-log canonical (slc) pairs and we refer the reader to \cite[Section 5]{Kol13} for their basic properties.

\begin{thm}\label{thm: special abundance}
Let $(X,B)$ be a $\mathbb Q$-factorial projective lc pair satisfying the following.
    \begin{enumerate}
        \item $X$ is klt and $B$ is reduced.
        \item $K_X+B$ is nef and $K_X+B\sim_{\mathbb R}D\geq 0$ for some $D$ such that $\Supp D\subset\Supp B$.
        \item $\kappa_{\sigma}(K_X+B)=1$ and $\kappa(K_X+B)\geq 0$.
        \item There exists an irreducible component $S$ of $D$ such that $S$ does not intersect $B-S$.
    \end{enumerate}
    Then $\kappa(K_X+B)=1$.
\end{thm}
\begin{proof}  
Note that by (2) and (4), we have
$$\kappa(K_X+B)=\kappa(D)\geq \kappa(S).$$
Hence to prove the assertion, it suffices to show that $\kappa(S)\geq 1$.
By \cite[Page 27]{BPV84} we can take $U$ to be an open subset of $X$ which retracts to $S$. 
Let $(S,B_S)$ be the pair induced by adjunction
$$K_S+B_S:=(K_U+B)|_S= (K_U+S)|_S= K_S+ \mathrm{Diff}_S(0),$$
where the second equality holds by (4), and $\mathrm{Diff}_S(0)$ is the different $\mathbb Q$-divisor. 
By \cite[16.9.1]{Kol+92}, $(S,B_S)$ is slc. Moreover, by Lemma \ref{lem: restrictio of nu=1 is nu=0}, $\kappa_{\sigma}(K_S+B_S)=0$ and $K_S+B_S$ is nef, so $K_S+B_S\equiv 0$. Hence by \cite[Theorem 1.5]{Gon13}, $K_S+B_S\sim_{\mathbb Q}0$.
Note that there exist two positive integers $a$ and $b$ such that $aK_U\sim bS$ and $b>-a$. Let $c:=\gcd(a,b)$. Then there exist integers $a'$, $b'$, $a_1$ and $b_1$ such that $a=ca'$, $b=cb'$ and $aa_1+bb_1=c$. Let $G=a_1 S+b_1 K_U$. Then we have 
$$c(S-a'G)=(aa_1+bb_1)S-a(a_1S+b_1K_U)=b_1(bS-aK_U)\sim 0,$$
$$c(K_U-b'G)=(aa_1+bb_1)K_U-b(a_1S+b_1K_U)=a_1(aK_U-bS)\sim 0,$$
and
$$c(K_U+S-(a'+b')G)=c(S-a'G)+c(K_U-b'G)\sim 0.$$
In particular, $c(S-a'G)$, $c(K_U-b'G)$ and $c(K_U+S-(a'+b')G)$ are Cartier divisors with first Chern classes equal to $0$ in $H^2(U,\mathbb Z)$. Moreover, since $(K_U+S)|_S= K_S+B_S\sim_{\mathbb Q} 0$, there exists a positive integer $m$ such that $mc(K_U+S)$ is Cartier, whose first Chern class is equal to $0$ in $H^2(U,\mathbb Z)$ by the isomorphism $H^2(U,\mathbb Z)\cong H^2(S,\mathbb Z)$. Hence $mc(a'+b')G\sim mc(K_U+S)$ is also Cartier, and its first Chern class is equal to $0$ in $H^2(U,\mathbb Z)$.
Applying \cite[11.3.6 Lemma]{Kol+92} to $S-a'G$, $K_U-b'G$ and $G$, we get
a finite cover $\pi: \tilde U\rightarrow U$, \'etale in codimension $1$, 
such that if we write $K_{\tilde U}:=\pi^*K_U$, $\tilde S:=\pi^*S$ and $\tilde G:=\pi^\ast G$,  the following conditions hold:
\begin{itemize}
    \item[(i)] $\tilde G$ is Cartier, and $\tilde G|_S\sim 0$.
    \item[(ii)] $\tilde S\sim a'\tilde G$ and $K_{\tilde U}\sim b'\tilde G$.
    \item[(iii)]  $(\tilde U,\tilde S)$ is a lc pair with $\tilde U$ klt and $\tilde S$ proper and reduced.
    \item[(iv)] $\omega_{\tilde S}\sim \mathcal O_{\tilde U}(K_{\tilde U}+\tilde S)|_{\tilde S}\sim \mathcal{O}_{\tilde S}$.
\end{itemize}
Here (iii) holds by \cite[Proposition 5.20]{KM98} since $\pi$ is \'etale in codimension $1$, and (iv) holds by (i) and (ii).

 We will prove that $\tilde S$ moves infinitesimally in $\Tilde U$. To this end, we use \cite[11.3.7 Lemma]{Kol+92} to conclude this. Note that since by (iii) $\tilde U$ is klt, $\tilde U$ is Cohen-Macaulay. To apply \cite[11.3.7 Lemma]{Kol+92}, it suffices to verify that the restriction
$$H^p(\tilde S_n, \mathcal{O}_{\tilde S_n}) \rightarrow H^p(\tilde S, \mathcal{O}_{\tilde S})$$
is surjective for every $p\in\mathbb N$, where $\tilde S_n$ are the analytic subspace of $\tilde U$ defined by the sheaf of ideals $\mathcal{O}_{\tilde U}(-n\tilde S)$. 
 By considering the commutative diagram 
 $$
\xymatrix{
H^p(\tilde S_n, \mathbb C) \ar[d] \ar[r]  & H^p(\tilde S_n, \mathcal{O}_{\tilde S_n}) \ar[d]
\\
H^p(\tilde S, \mathbb C) \ar[r] & H^p(\tilde S, \mathcal{O}_{\tilde S})
}
$$
and the isomorphism $H^p(\tilde S_n,  \mathbb C) \simeq H^p(\tilde S,  \mathbb C)$, it is sufficient to prove that
$$
H^p(\tilde S,  \mathbb C) \rightarrow H^p(\tilde S, \mathcal{O}_{\tilde S})
$$ 
is surjective. By (iii) and \cite[16.9.1]{Kol+92}, $\tilde S$ is slc. 
By \cite[Corollary 6.32]{Kol13} (see also \cite[Theorem 1.4]{KK10}), $\tilde S$ is Du Bois, and so $
H^p(\tilde S,  \mathbb C) \rightarrow H^p(\tilde S, \mathcal{O}_{\tilde S})
$ is surjective (cf. \cite[Property 2.61.2]{Kol23}). Therefore, applying \cite[11.3.7 Lemma]{Kol+92}, we know that $\tilde S$ moves infinitesimally in $\Tilde U$. It implies that $\kappa(S)=\kappa(\tilde S)\geq 1,$
and we are done.
\end{proof}

\section{The projective case}

In this section, we prove Theorem \ref{thm: main sp tof} for projective lc pairs. More precisely, the main theorem of this section is the following:

\begin{thm}\label{thm: projective case}
 Assume Conjecture \ref{conj: stof with scaling weak} in dimension $\leq d$. Let $(X,B)$ be a projective lc pair of dimension $d$ such that $\kappa_{\iota}(K_X+B)\geq 0$ and $\kappa_{\sigma}(K_X+B)\leq 1$. Then $(X,B)$ has a good minimal model. In particular, $\kappa_{\iota}(K_X+B)=\kappa_{\sigma}(K_X+B)$.
\end{thm}

The proof proceeds in several steps: We first reduce to the case of a projective lc pair $(X,B)$ such that $K_X + B \sim_{\mathbb{R}} D \geq 0$ with $\Supp D \subset \Supp \lfloor B \rfloor$, while preserving the numerical dimension and Kodaira dimension of the original pair (Lemmas \ref{lem: klt to special log smooth} and \ref{lem: klt to nef dlt}). Then we further reduce to the case when $\Supp D=\Supp\lfloor B\rfloor$ (Lemma \ref{lem: d subset b to d=b}), to the case where $B$ is reduced (Lemmas \ref{lem: reduction to lower dimension if not pseudo-effective} and \ref{lem: reduce to reduced boundary}), and finally to the case where $B$ has an isolated irreducible component (Lemma \ref{lem: reduce to disconnect s}) and apply the special case established earlier.

\begin{lem}\label{lem: klt to special log smooth}
Let $(X,B)$ be a projective lc pair of dimension $d$ such that $\kappa_{\iota}(K_X+B)\geq 0$. Then there exists a log smooth pair $(X',B')$ of dimension $d$ satisfying the following.
\begin{enumerate}
\item $K_{X'}+B'\sim_{\mathbb R}D'\geq 0$ for some $D'$ such that $\Supp D'\subset\Supp\lfloor B'\rfloor$.
\item $\kappa_{\sigma}(K_X+B)=\kappa_{\sigma}(K_{X'}+B')$ and $\kappa_{\iota}(K_X+B)=\kappa_{\iota}(K_{X'}+B')$.
\end{enumerate}
\end{lem}
\begin{proof}
 Write $K_X+B\sim_{\mathbb R}D\geq 0$. We let $h: X'\rightarrow X$ be a log resolution of $(X,\Supp B\cup\Supp D)$ and let $$B':=\Exc(h)+h^{-1}_*\Supp D+h^{-1}_*(B-B\wedge\Supp D).$$
Write $K_{X'}+\bar B':=h^*(K_{X'}+B')$. Then
$$K_{X'}+B'=K_{X'}+\bar B'+((h^{-1}_*B+\Exc(h))-\bar B')+h^{-1}_*(\Supp D-B\wedge\Supp D).$$
We have
$$L_1:=(h^{-1}_*B+\Exc(h))-\bar B'\geq 0\text{ and }\Supp L_1\subset\Exc(h),$$
and
$$L_2:=h^{-1}_*(\Supp D-B\wedge\Supp D)\geq 0 \text{ and } \Supp L_2\subset h^{-1}_*\Supp D.$$
We let $D':=h^*D+L_1+L_2$. Since
$$\Supp h^*D\subset h^{-1}_*\Supp D\cup\Exc(h),$$
we have
$$\Supp D'\subset h^{-1}_*\Supp D\cup\Exc(h)\subset\Supp\lfloor B'\rfloor.$$
This implies (1).

(3) Since $(K_{X'}+B')-h^*(K_X+B)\geq 0$, we have  
$$\kappa_{\sigma}(K_X+B)\leq\kappa_{\sigma}(K_{X'}+B') \text{ and }\kappa_{\iota}(K_X+B)\leq\kappa_{\iota}(K_{X'}+B').$$
Let $s\gg 0$ be a real number such that $sD\geq\Supp D$ and $\bar B'+s\Exc(h)\geq h^{-1}_*B+\Exc(h)$. Let
$$K_{X'}+\tilde B':=h^*(K_{X}+B+sD)+s\Exc(h),$$
then 
$$\tilde B'=\bar B+s\Exc(h)+sh^*D\geq h^{-1}_*B+\Exc(h)+h^{-1}_*\Supp D\geq B',$$
hence
$$\kappa_{\sigma}(K_{X'}+\tilde B')\geq\kappa_{\sigma}(K_{X'}+B') \text{ and }\kappa_{\iota}(K_{X'}+\tilde B')\geq\kappa_{\iota}(K_{X'}+B').$$
However, by Lemma \ref{lem: property of numerical and Iitaka dimension}, we have
$$\kappa_{\sigma}(K_{X'}+\tilde B')=\kappa_{\sigma}(K_{X}+B+sD)=\kappa_{\sigma}((1+s)(K_X+B))=\kappa_{\sigma}(K_X+B)$$
and
$$\kappa_{\iota}(K_{X'}+\tilde B')=\kappa_{\iota}(K_{X}+B+sD)=\kappa_{\iota}((1+s)(K_X+B))=\kappa_{\iota}(K_X+B).$$
Thus (3) follows.
\end{proof}

\begin{lem}\label{lem: klt to nef dlt}
 Assume Conjecture \ref{conj: stof with scaling weak} dimension $d$. Let $(X,B)$ be a projective lc pair of dimension $d$ such that $\kappa_{\iota}(K_X+B)\geq 0$. Then there exists a projective $\mathbb Q$-factorial dlt pair $(X'',B'')$ of dimension $d$ satisfying the following.
\begin{enumerate}
    \item $K_{X''}+B''\sim_{\mathbb R}D''\geq 0$ for some $D''$.
    \item $\Supp D''\subset\Supp\lfloor B''\rfloor$.
    \item $\kappa_{\sigma}(K_X+B)=\kappa_{\sigma}(K_{X''}+B'')$ and $\kappa_{\iota}(K_X+B)=\kappa_{\iota}(K_{X''}+B'')$.
    \item $K_{X''}+B''$ is nef.
\end{enumerate}
\end{lem}
\begin{proof}
By Lemma \ref{lem: klt to special log smooth} there exists a projective log smooth pair $(X',B')$ of dimension $d$ satisfying the following:
\begin{itemize}
    \item $K_{X'}+B'\sim_{\mathbb R}D'\geq 0$ for some $D'$.
    \item $\Supp D'\subset\Supp\lfloor B'\rfloor$.
    \item $\kappa_{\sigma}(K_X+B)=\kappa_{\sigma}(K_{X'}+B')$ and $\kappa_{\iota}(K_X+B)=\kappa_{\iota}(K_{X'}+B')$.
\end{itemize}
By our assumption, we may run a $(K_{X'}+B')$-MMP so that after finitely many steps we obtain a birational map 
$$\phi: (X',B')\dashrightarrow (X'',B'')$$
such that the MMP terminates near $\lfloor B''\rfloor$. However, this MMP is also a $D'$-MMP. Since  $\Supp D'\subset\Supp\lfloor B'\rfloor$, this MMP terminates at $(X'',B'')$. In particular, $K_{X''}+B''$ is nef. Let $D'':=\phi_*D$. By Lemma \ref{lem: property of numerical and Iitaka dimension}, $(X'',B'')$ and $D''$ satisfy our requirements.
\end{proof}

\begin{lem}\label{lem: d subset b to d=b}
 Assume Conjecture \ref{conj: stof with scaling weak} dimension $d$ and Theorem \ref{thm: projective case} holds in dimension $d-1$. Let $(X,B)$ be a projective $\mathbb Q$-factorial dlt pair of dimension $d$ satisfying the following.
\begin{enumerate}
    \item $K_{X}+B\sim_{\mathbb R}D\geq 0$ for some $D$.
    \item $\Supp D\subset\Supp\lfloor B\rfloor$.
    \item $\kappa_{\sigma}(K_X+B)=1$ and $\kappa_{\iota}(K_X+B)=0$.
    \item $K_X+B$ is nef.
\end{enumerate}
Then $\Supp D=\Supp B$.
\end{lem}
\begin{proof}
Suppose that the lemma does not hold. Then there exists an irreducible component $S$ of $\lfloor B\rfloor$ that is not an irreducible component of $D$. There are two cases.

\medskip

\noindent\textbf{Case 1.} $K_X+B-\epsilon S$ is pseudo-effective for some rational number $\epsilon>0$. Since $\kappa_{\sigma}(K_X+B)=1$ and
$$K_X+B-\frac{\epsilon}{2}S=(K_X+B-\epsilon S)+\frac{\epsilon}{2}S,$$
By Theorem \ref{thm: lp18 r coefficient}, there exists $E\geq 0$ such that
$$K_X+B-\frac{\epsilon}{2}S\equiv E.$$
By Lemma \ref{lem: ckp12 real}, there exists $D'\geq 0$ on $X$ such that $K_X+B-\frac{\epsilon}{2}S\sim_{\mathbb R}D'$. Since $\kappa_{\iota}(K_X+B)=0$ and $K_X+B\sim_{\mathbb R}D\geq 0$, $D$ is the only element in $|K_X+B|_{\mathbb R}$. However, $$0\leq D'+\frac{\epsilon}{2}S\sim_{\mathbb R}K_X+B,$$
this is not possible as $S$ is not an irreducible component of $D$.

\medskip

\noindent\textbf{Case 2.}  $K_X+B-\epsilon S$ is not pseudo-effective for any rational number $\epsilon>0$. By Lemma \ref{lem: hl22 3.21}, we may pick a rational number $0<\epsilon\ll 1$ and run a $(K_X+B-\epsilon S)$-MMP with scaling of an ample divisor which terminates with a Mori fiber space $\pi: X'\rightarrow Z$ with induced birational map $\phi: X\dashrightarrow X'$, and this MMP is $(K_X+B)$-trivial. Let $B':=\phi_*B$ and $\tilde S':=\phi_*S$. Then this MMP is also a sequence of steps of a $(-S)$-MMP, so $S$ is not contracted by this MMP and $\tilde S'$ is ample$/Z$, hence $\tilde S'$ is horizontal$/Z$. Let $S'$ be the normalization of $\tilde S'$, $\pi_S: S'\rightarrow Z$ the induced projective surjective morphism, and let
$$K_{S'}+B_{S'}:=(K_{X'}+B')|_{S'}.$$
Then $(S',B_{S'})$ is a projective lc pair of dimension $d-1$ and there exists an $\mathbb R$-divisor $L$ on $Z$ such that
$$K_{X'}+B'=\pi^*L\text{ and }K_{S'}+B_{S'}=\pi_S^*L.$$
Thus
$$\kappa_{\sigma}(K_{S'}+B_{S'})=\kappa_{\sigma}(L)=\kappa_{\sigma}(K_{X'}+B')=1$$
and
$$\kappa_{\iota}(K_{S'}+B_{S'})=\kappa_{\iota}(L)=\kappa_{\iota}(K_{X'}+B')=0$$
which contradicts Theorem \ref{thm: projective case} in dimension $d-1$.
\end{proof}

\begin{lem}\label{lem: reduction to lower dimension if not pseudo-effective}
    Let $(X,B)$ be a projective $\mathbb Q$-factorial dlt pair of dimension $d$ satisfying the following.
    \begin{itemize}
        \item $K_X+B\sim_{\mathbb R}D\geq 0$ for some $D$ such that $\Supp D=\Supp\lfloor B\rfloor$. 
        \item $K_X+B$ is nef.
        \item $K_X+B-\epsilon\{B\}$ is not pseudo-effective for any $\epsilon>0$.
    \end{itemize}
Then either $K_X+B\sim_{\mathbb R}0$, or there exists a klt pair $(Z,B_Z)$ satisfying the following.
\begin{enumerate}
    \item $0<\dim Z<\dim X$.
    \item $\kappa_{\sigma}(K_X+B)=\kappa_{\sigma}(K_Z+B_Z)$ and $\kappa_{\iota}(K_X+B)=\kappa_{\iota}(K_Z+B_Z)$.
\end{enumerate}
\end{lem}
\begin{proof}
Pick $0<\epsilon\ll 1$ and run a $(K_X+B-\epsilon\{B\})$-MMP with scaling of an ample divisor which terminates with a Mori fiber space $f: X'\rightarrow Z$. Let $\phi: (X,B)\dashrightarrow (X',B')$ be this MMP. Since $0<\epsilon\ll 1$, by Lemma \ref{lem: hl22 3.21}, $\phi$ is $(K_X+B)$-trivial. In particular, $K_{X'}+B'\sim_{\mathbb R,Z}0$. If $\dim Z=0$ then $K_{X'}+B'\sim_{\mathbb R}0$, and since $\phi$ is $(K_X+B)$-trivial, we have $K_X+B\sim_{\mathbb R}0$ and we are done. Thus may assume that $\dim Z>0$. It is clear that $\dim Z<\dim Z$. This implies (1).

To prove (2), let $\Delta:=K_X+B-\delta D$ for some $0<\delta\ll 1$. Since $(X,B)$ is $\mathbb Q$-factorial dlt and $\Supp D=\Supp\lfloor B\rfloor$, $(X,\Delta)$ is klt. Let $\Delta'$ be the image of $\Delta$ on $X'$. Since $K_X+\Delta\sim_{\mathbb R}(1-\delta)(K_X+B)$, this MMP is $(K_X+\Delta)$-trivial, hence $(X',\Delta')$ is klt. By Lemma \ref{lem: ambro05 r coefficient} there exists a klt pair $(Z,B_Z)$ such that
$$(1-\delta)(K_{X'}+B')\sim_{\mathbb R}K_{X'}+\Delta'\sim_{\mathbb R}f^*(K_Z+B_Z).$$
(2) follows.
\end{proof}

\begin{lem}\label{lem: reduce to reduced boundary}
 Assume Conjecture \ref{conj: stof with scaling weak} in dimension $d$. Let $(X,B)$ be a projective $\mathbb Q$-factorial dlt pair of dimension $d$ satisfying the following.
    \begin{itemize}
        \item $K_X+B\sim_{\mathbb R}D\geq 0$ for some $D$ such that $\Supp D=\Supp\lfloor B\rfloor$. 
        \item $\kappa_{\sigma}(K_X+B)=1$.
        \item $K_X+B$ is nef.
        \item $K_X+B-\epsilon\{B\}$ is pseudo-effective for some $\epsilon>0$.
    \end{itemize}
Then either $(X,B)$ has a good minimal model, or there exists a $\mathbb Q$-factorial projective dlt pair $(X'',B'')$ satisfying the following.
\begin{enumerate}
    \item $K_{X''}+B''\sim_{\mathbb R}D''\geq 0$ for some $D''$ such that $\Supp D''=\Supp B''$.
    \item $\kappa_{\sigma}(K_{X''}+B'')=1$ and $\kappa(K_{X''}+B'')=0$.
    \item $K_{X''}+B''$ is nef.
    \item $B''$ is reduced.
\end{enumerate}
\end{lem}
\begin{proof}
The proof is a little bit long so we separate it into several steps.

\medskip

\noindent\textbf{Step 1.} In this step we introduce $\kappa_0:=\kappa_{\iota}(K_X+B+P)$ and show that it is a constant in $\{0,1\}$ for any $P\geq 0$ such that $\Supp P= \Supp\{ B\}$.

Since $K_X+B-\epsilon\{B\}$ is pseudo-effective, by Lemma \ref{lem: compare kappa}, we have 
$$\kappa_{\sigma}(K_X+B+t\{B\})=\kappa_{\sigma}(K_X+B)=1,$$
so
$$1\geq\kappa_{\iota}(K_X+B+t\{B\})\geq\kappa_{\iota}(K_X+B)\geq 0.$$
By Lemma \ref{lem: compare kappa}, 
for any $t>0$, we have that $\kappa_0:=\kappa_{\iota}(K_X+B+t\{B\})$ is a constant in $\{0,1\}$ for any $t>0$. Moreover, for any $\mathbb R$-divisor $P\geq 0$ such that $\Supp P=\Supp\{B\}$, there exist real numbers $t_1>t_2>0$ such that $t_1\{B\}\geq P\geq t_2\{B\}$, so by  Lemma \ref{lem: compare kappa} again, we have $\kappa_{\sigma}(K_X+B+P)=1$ and $\kappa_{\iota}(K_X+B+P)=\kappa_0$.

\medskip

\noindent\textbf{Step 2.} In this step we construct a projective log smooth pair $(X',B')$ with reduced boundary which satisfies certain properties.

Since $(X,B)$ is $\mathbb Q$-factorial dlt, there exists a log resolution $h: X'\rightarrow X$ of $(X,B)$ which only extract non-lc places of $(X,B)$. Let $B':=h^{-1}_*\Supp B+\Exc(h)$. Since $(X,B)$ is dlt, we have
$$K_{X'}+B'\geq h^*(K_X+B+t\{B\})$$
for some $0<t\ll 1$, hence
$$\kappa_{\sigma}(K_{X'}+B')\geq\kappa_{\sigma}(K_X+B+t\{B\})=1$$
and
$$\kappa_{\iota}(K_{X'}+B')\geq\kappa_{\iota}(K_X+B+t\{B\})=\kappa_0.$$
Let $s\gg 0$ be a positive real number such that $B+s\{B\}\geq\Supp B$ and $h^*(K_X+B)+s\Exc(h)\geq K_{X'}+h^{-1}_*B+\Exc(h)$, then
$$h^*(K_X+B+s\{B\})+s\Exc(h)\geq K_{X}+h^{-1}_*\Supp B$$
and
$$h^*(K_X+B+s\{B\})+s\Exc(h)\geq h^*(K_X+B)+s\Exc(h)\geq K_{X'}+h^{-1}_*B+\Exc(h),$$
hence 
$$h^*(K_X+B+s\{B\})+s\Exc(h)\geq K_X+B'.$$
By Lemma \ref{lem: property of numerical and Iitaka dimension} and since $B'$ is a Weil divisor, we have
$$\kappa_{\sigma}(K_{X'}+B')\leq\kappa_{\sigma}(h^*(K_X+B+s\{B\})+s\Exc(h))=\kappa_{\sigma}(K_X+B+s\{B\})=1$$
and
$$\kappa(K_{X'}+B')=\kappa_{\iota}(K_{X'}+B')\leq\kappa_{\iota}(h^*(K_X+B+s\{B\})+s\Exc(h))=\kappa_{\iota}(K_X+B+s\{B\})=\kappa_0.$$
Therefore, $\kappa_{\sigma}(K_{X'}+B')=1$ and $\kappa(K_{X'}+B')=\kappa_0$.

Moreover, write
$$K_{X'}+h^{-1}_*B+E:=h^*(K_X+B),$$
then $\lfloor E\rfloor\leq 0$ and $E$ is exceptional$/X$. Therefore, 
$$K_{X'}+B'\sim_{\mathbb R}h^*D+(\Exc(h)-E)+h^{-1}_*(\Supp B-B)=:D'\geq 0.$$
It is clear that $\Supp D'\subset\Supp B'$. Moreover, since $\Supp D=\lfloor B\rfloor$, $\Supp(\Supp B-B)=\Supp\{B\}$, and $h_*D'=D+\Supp B-B$, we have $\Supp h_*D'=\Supp B$. Since $\Supp(\Exc(h)-E)=\Exc(h)$, we have $\Supp D'=\Supp B'$.

\medskip

\noindent\textbf{Step 3.} In this step we construct $(X'',B'')$ and $D''$ and show that they satisfy all the required properties except $\kappa(K_{X''}+B'')=0$.

By our assumption, we may run a $(K_{X'}+B')$-MMP so that after finitely many steps we obtain a birational map 
$$\phi: (X',B')\dashrightarrow (X'',B'')$$
such that the MMP terminates near $\lfloor B''\rfloor$. However, this MMP is also a $D'$-MMP. Since  $\Supp D'\subset\Supp\lfloor B'\rfloor$, this MMP terminates at $(X'',B'')$. In particular, $K_{X''}+B''$ is nef. Let $D'':=\phi_*D$. Then $\Supp D''=\Supp B''$.

By our construction, $(X'',B'')$ and $D''$ satisfy all the required properties except $\kappa(K_{X''}+B'')=0$, and by Lemma \ref{lem: property of numerical and Iitaka dimension}, we have $\kappa(K_{X''}+B'')=\kappa_0$. 

\medskip

\noindent\textbf{Step 4.} In this step we conclude the proof. We may assume that $\kappa_0=1$.

We pick $0<t\ll 1$ and $0<\delta\ll 1$ such that $(X,B+\frac{t}{1-\epsilon}\{B\})$ is dlt and $\lfloor B+\frac{t}{1-\epsilon}\{B\}\rfloor=\lfloor B\rfloor$. In particular, $(X,B-\delta D+t\{B\})$ is klt. Since 
$$K_X+B-\delta D+t\{B\}\sim_{\mathbb R}(1-\delta)\left(K_X+B+\frac{t}{1-\delta}\{B\}\right),$$
we have
$$\kappa_{\sigma}(K_X+B-\delta D+t\{B\})=\kappa_{\iota}(K_X+B-\delta D+t\{B\})=1.$$
Since $(X,B-\delta D+t\{B\})$ is klt, by Lemma \ref{lem: hh20 2.13}, $(X,B-\delta D+t\{B\})$ has a good minimal model, hence $(X,B+\frac{t}{1-\epsilon}\{B\})$ has a good minimal model. By Theorem \ref{thm: mz23 1.4}, $(X,B)$ has a good minimal model and we are done. 
\end{proof}

\begin{lem}\label{lem: reduce to disconnect s}
 Assume Conjecture \ref{conj: stof with scaling weak} in dimension $d$. Let $(X,B)$ be a projective $\mathbb Q$-factorial dlt pair of dimension $d$ satisfying the following.
    \begin{itemize}
    \item $X$ is klt and $B$ is reduced.
    \item $K_X+B\sim_{\mathbb R}D\geq 0$ for some $D$ such that $\Supp D=\Supp B$. 
    \item $\kappa_{\sigma}(K_X+B)=1$ and  $\kappa(K_X+B)=0$. 
    \item $K_X+B$ is nef.    
    \end{itemize}
    Then there exists a $\mathbb Q$-factorial lc pair $(X',B')$ crepant to $(X,B)$ satisfying the following.
    \begin{enumerate}
    \item $X'$ is klt and $B'$ is reduced.
    \item $K_{X'}+B'\sim_{\mathbb R}D'\geq 0$ for some $D'$ such that $\Supp D'=\Supp B'$. 
   \item  There exists an irreducible component $S'$ of $B'$ such that $S'$ does not intersect $B'-S'$.
    \end{enumerate}
    In particular, $\kappa_{\sigma}(K_{X'}+B')=1$ and  $\kappa(K_{X'}+B')=0$.
\end{lem}
\begin{proof}
Since $\kappa_{\sigma}(K_X+B)=1$, we have $D\not=0$. Let $S$ be any irreducible component of $D$. Thus we may write $D=aS+D_0$ for some $a>0$, where $\Supp D_0=\Supp(B-S)$. By Lemma \ref{lem: hl22 3.21}, we may pick $0<\epsilon\ll 1$ such that any sequence of steps of a $(K_X+B-\epsilon S)$-MMP is $(K_X+B)$-trivial. Since $\kappa_{\iota}(K_X+B-\epsilon S)=\kappa((a-\epsilon)S+D_0)\geq 0$, by our assumption, we may run a $(K_X+B-\epsilon S)$-MMP which, after finitely many steps, we obtain a birational map 
$$\phi: (X',B')\dashrightarrow (X',B'),$$
such that the MMP terminates near $\lfloor B'-S'\rfloor$, where $S':=\phi_*S$. Since this MMP is a $(-S)$-MMP and is $(aS+D_0)$-trivial, it is also a $D_0$-MMP. Since $\Supp D_0=\Supp(B-S)$, this MMP terminates at $(X',B')$. In particular, $K_{X'}+B'$ is nef. Let $D':=\phi_*D$.

We show that $(X',B')$ satisfies our requirements. 

(1) By our construction, $\phi$ is $(K_X+B)$-trivial, hence $(X',B')$ is crepant to $(X,B)$. In particular, $(X',B')$ is lc. Since $(X,B-\epsilon S)$ is $\mathbb Q$-factorial dlt, $(X',B'-\epsilon S')$ is $\mathbb Q$-factorial dlt. Thus $X'$ is klt. Since $B$ is reduced, $B'$ is reduced. 

(2) It is immediate from our construction.

(3) Suppose not, then there exists an irreducble componet $\tilde S'$ of $B'$ which intersects $S'$.

By Theorem \ref{thm: abundance dim 3}, we have $d\geq 4$. Let $H$ be a very ample divisor on $X'$ and let $H_1,\dots,H_{d-2}\in |H|$ be general elements. We let $V=S'\cap\tilde S'$, then $\dim V=\dim X-2$. Let $C:=H_1\cap H_2\dots\cap H_{d-3}\cap V$, then $C$ is curve (not necessarily irreducible) contained in $S'\cap\tilde S'$. Let $\tilde C:=H_1\cap H_2\dots\cap H_{d-2}\cap\tilde S'$, then $\tilde C$ is a curve contained in $\tilde S'$. We have
$$0<H\cdot C=S'\cdot\tilde C.$$
By Lemma \ref{lem: restrictio of nu=1 is nu=0}, $(K_{X'}+B')|_{\tilde S'}\equiv 0$, hence we have 
$$(K_{X'}+B')\cdot\tilde C=(K_{X'}+B')|_{\tilde S'}\cdot\tilde C=0$$
as $\tilde C\subset\tilde S'$. Therefore, $(K_{X'}+B'-\delta S')\cdot\tilde C<0$ for any $\delta>0$. This is not possible since $K_{X'}+B'-\epsilon S'$ is nef.
\end{proof}

\begin{proof}[Proof of Theorem \ref{thm: projective case}]
We may apply induction on dimension and assume that the theorem holds in dimension $<d$. By Theorem \ref{thm: num 0 abundance}, we may assume that $\kappa_{\iota}(K_X+B)=0$ and $\kappa_{\sigma}(K_X+B)=1$. By Lemma \ref{lem: klt to nef dlt}, there exists a projective $\mathbb Q$-factorial dlt pair $(X_1,B_1)$ of dimension $d$ satisfying the following.
\begin{itemize}
    \item $K_{X_1}+B_1\sim_{\mathbb R}D_1\geq 0$ for some $D_1$ such that $\Supp D_1\subset\Supp\lfloor B_1\rfloor$.
    \item $\kappa_{\sigma}(K_{X_1}+B_1)=1$ and $\kappa_{\iota}(K_{X_1}+B_1)=0$.
    \item $K_{X_1}+B_1$ is nef.
\end{itemize}
By Lemma \ref{lem: d subset b to d=b} and induction hypothesis, we have $\Supp D_1=\Supp\lfloor B_1\rfloor$. By Lemma \ref{lem: reduction to lower dimension if not pseudo-effective} and induction hypothesis, we may assume that $K_{X_1}+B_1-\epsilon\{B_1\}$ is pseudo-effective for some $\epsilon>0$. By Lemma \ref{lem: reduce to reduced boundary}, we may assume that there exists a $\mathbb Q$-factorial dlt pair $(X_2,B_2)$ satisfying the following.
\begin{itemize}
    \item $K_{X_2}+B_2\sim_{\mathbb R}D_2\geq 0$ for some $D_2$ such that $\Supp D_2=\Supp B_2$.
    \item $\kappa_{\sigma}(K_{X_2}+B_2)=1$ and $\kappa(K_{X_2}+B_2)=0$.
    \item $K_{X_2}+B_2$ is nef.
    \item $B_2$ is reduced.
\end{itemize}

By Lemma \ref{lem: reduce to disconnect s}, we may further assume that there exists an irreducible component $S_2$ of $B_2$ such that $S_2$ does not intersect $B_2-S_2$. This contradicts Theorem \ref{thm: special abundance}.
\end{proof}

\section{Proof of the main theorems}\label{sec: proof of main theorems}

\begin{proof}[Proof of Theorem \ref{thm: main sp tof}]
    It is a special case of Theorem \ref{thm: projective case}.
\end{proof}

\begin{thm}\label{thm: lc abundance dim 5 nu 1}
   For any projective lc pair $(X,B)$ of dimension $\leq 5$ such that $\kappa_{\iota}(K_X+B)\geq 0$ and $\kappa_{\sigma}(K_X+B)\leq 1$, $(X,B)$ has a good minimal model. 
\end{thm}
\begin{proof}
Conjecture \ref{conj: stof with scaling weak} holds in dimension $\leq 5$ by \cite[Lemmas 3.6, 3.8]{Bir10}. So the theorem follows from Theorem \ref{thm: projective case}.
\end{proof}

\begin{proof}[Proof of Theorem \ref{thm: main dim 5}]
It is a special case of Theorem \ref{thm: lc abundance dim 5 nu 1}.
\end{proof}

\begin{cor}\label{cor: lp18 klt}
Let $(X,B)$ be a klt pair of dimension $\leq 4$ such that $\chi(\mathcal{O}_X)\not=0$ and $0\leq\kappa_{\sigma}(K_X+B)\leq 1$. Then $(X,B)$ has a good minimal model. In particular, $\kappa_{\sigma}(K_X+B)=\kappa_{\iota}(K_X+B)$.
\end{cor}
\begin{proof}
By Theorem \ref{thm: abundance dim 3}, we may assume that $\dim X=4$. By Theorem \ref{thm: lc abundance dim 5 nu 1}, we only need to show that $\kappa_{\iota}(K_X+B)\geq 0$. Let $h: X'\rightarrow X$ be a log resolution of $(X,B)$ and $B':=h^{-1}_*B+(1-\epsilon)\Exc(h)$ for some $0<\epsilon\ll 1$. Since $X$ has rational singularities, $\chi(\mathcal{O}_{X'})=\chi(\mathcal{O}_X)\not=0$. Possibly replacing $(X,B)$ with $(X',B')$, we may assume that $(X,B)$ is log smooth. If $K_X$ is not pseudo-effective, then $X$ is uniruled by \cite{BDPP13} and we done by \cite[Corollary 1.2]{LM21}. Thus $K_X$ is pseudo-effective. We have
$$1=\kappa_{\sigma}(K_{X'}+B')\geq\nu(X)\geq 0.$$
By Theorem \ref{thm: num 0 abundance}, we may assume that $\nu(X)=1$. By \cite{Fuj04,Fuj05}, we may run a $K_X$-MMP which terminates with a minimal model $Y$ of $X$. Then $Y$ is terminal, $\chi(\mathcal{O}_Y)=\chi(\mathcal{O}_X)\not=0$, $\nu(Y)=\nu(X)=1$, and $\kappa(X)=\kappa(Y)$. By \cite[Theorem 6.7]{LP18}, 
$$\kappa_{\iota}(K_X+B)\geq\kappa(X)=\kappa(Y)\geq 0$$
and we are done.
\end{proof}

\begin{proof}[Proof of Corollary \ref{cor: lp18}]
  This is a special case of Corollary \ref{cor: lp18 klt}.  
\end{proof}

We are left to prove Theorem \ref{thm: main nv}. Here we prove a more general case of Theorem \ref{thm: main nv} by only assuming the non-vanishing conjecture in lower dimensions and for smooth projective varieties of numerical dimension one. We need the following lemma:

\begin{lem}\label{lem: non pe induces contraction}
Let $(X,B)$ be a projective lc pair and let $S$ be an irreducible component of $B$ such that $S$ is $\mathbb Q$-Cartier. Assume that $K_X+B$ is pseudo-effective and $K_X+B-\epsilon S$ is not pseudo-effective for any $\epsilon>0$. Then there exists a birational map $\varphi: X\dashrightarrow Y$, a pair $(Y,B_Y)$, and a contraction $f: Y\rightarrow Z$ satisfying the following.
\begin{enumerate}
    \item $S_Y:=\varphi_*S$ is an irreducible component of $B_Y$ with $\mult_{S_Y}B_Y=\mult_SB$.
    \item $K_Y+B_Y\sim_{\mathbb R,Z}0$.
    \item $\dim Y>\dim Z$ and $S_Y$ is horizontal$/Z$.
    \item $\kappa_{\sigma}(K_X+B)=\kappa_{\sigma}(K_Y+B_Y)$ and $\kappa_{\iota}(K_X+B)=\kappa_{\iota}(K_Y+B_Y)$.
    \item If $(X,B)$ is klt, then $(Y,B_Y)$ is $\mathbb Q$-factorial klt.
\end{enumerate}
\end{lem}
\begin{proof}
    By \cite[Theorem 1.7]{HH20}, we may run a $(K_X+B-\epsilon S)$-MMP$/U$ which terminates with a Mori fiber space$/U$: $f_{\epsilon}: X_{\epsilon}\rightarrow Z_{\epsilon}$. Let $\phi_{\epsilon}: X\dashrightarrow X_{\epsilon}$ be the induced birational map and let $B_{\epsilon}:=(\phi_{\epsilon})_*B,S_{\epsilon}:=(\phi_{\epsilon})_*S$. Since $K_X+B$ is pseudo-effective, $K_{X_{\epsilon}}+B_{\epsilon}$ is pseudo-effective for any $0<\epsilon\ll 1$. Since $\rho(X_{\epsilon}/Z_{\epsilon})=1$, for any $0<\epsilon\ll 1$, there exists $\delta_{\epsilon}\in [0,\epsilon)$ such that $K_{X_{\epsilon}}+B_{\epsilon}-\delta_{\epsilon}S_{\epsilon}\sim_{\Rr,Z_{\epsilon}}0$. By \cite[Theorem 1.1]{HMX14}, for any $0<\epsilon\ll 1$, $(X_{\epsilon},B_{\epsilon})$ is lc, hence $(X_{\epsilon},B_{\epsilon}-\delta_{\epsilon}S_{\epsilon})$ is lc. By \cite[Theorem 1.5]{HMX14}, for any $0<\epsilon\ll 1$, $\delta_{\epsilon}=0$ and $K_{X_{\epsilon}}+B_{\epsilon}\sim_{\Rr,Z_{\epsilon}}0$. In the following, we will fix $0<\epsilon\ll 1$ and denote $X':=X_{\epsilon}$, $Z':=Z_{\epsilon}$, $\phi:=\phi_{\epsilon}$, $B':=B_{\epsilon}$, and $S':=S_{\epsilon}$.

    Let $p: \bar X\rightarrow X$ and $q: \bar X\rightarrow X'$ be a resolution of indeterminacy of $\phi$ such that $p$ is a log resolution of $(X,\Supp B)$ and $q$ is a log resolution of $(X',\Supp B')$. Then we may write
$$K_{\bar X}+\bar B=p^*(K_X+B)+E$$
for some $p$-exceptional $\Rr$-divisor $E\geq 0$ and dlt pair $(\bar X,\bar B)$ such that $\bar B\wedge E=0$, and write
$$K_{\bar X}+\bar B=q^*(K_{X'}+B')+F-G$$
for some $q$-exceptional $\Rr$-divisors $F\geq 0,G\geq 0$ such that $F\wedge G=0$. By \cite[Theorem 3.5]{Bir12}, we may run a $(K_{\bar X}+\bar B)$-MMP$/X'$ which terminates with a model $\bar X'$ with induced birational map $\psi: \bar X\dashrightarrow\bar X'$. It is clear that $\psi$ only contracts $q$-exceptional prime divisors, and in particular, $p^{-1}_*S$ is not contracted by $\psi$. Let $\bar B',\bar G'$be the images of $\bar B, G$ on $\bar X'$ respectively and $q':\bar X'\rightarrow X'$ the induced morphism. Then $(\bar X',\bar B')$ is $\mathbb Q$-factorial dlt,
$$\kappa_{\sigma}(K_{\bar X'}+\bar B')=\kappa_{\sigma}(K_{\bar X}+\bar B)=\kappa_{\sigma}(K_X+B),$$
and
$$\kappa_{\iota}(K_{\bar X'}+\bar B')=\kappa_{\iota}(K_{\bar X}+\bar B)=\kappa_{\iota}(K_X+B).$$
In particular, $K_{\bar X'}+\bar B'$ is pseudo-effective. Since
$$K_{\bar X'}+\bar B'+\bar G'=q'^*(K_{X'}+B')\sim_{\Rr,Z'}0,$$
$\bar G'=0$ over the generic point of $Z'$, and by \cite[Theorem 1.1]{Has19a}, we may run a $(K_{\bar X'}+\bar B')$-MMP$/Z'$ with scaling of an ample divisor which terminates with a good minimal model$/Z'$ $(Y,B_Y)/Z'$ of $(\bar X',\bar B')/Z'$. This MMP is an isomorphism over the generic point of $Z'$, and in particular, $q^{-1}_*S'$ is not contracted by this MMP. Let $f: Y\rightarrow Z$ be the ample model$/Z'$ of $K_Y+B_Y$, $S_Y$ the strict transform of $S$ on $Y$, and $\varphi: X\dashrightarrow Y$ the induced birational map.

We show that $\varphi,(Y,B_Y)$ and $f$ satisfy our requirements. We have already checked that $\varphi$ does not extract $S$. Since $p^{-1}_*S$ is an irreducible component of $\bar B$ and $B_Y$ is the image of $\bar B$ on $Y$, $S_Y$ is an irreducible component of $B_Y$ and $\mult_{S_Y}B_Y=\mult_SB$. This implies (1).

(2) is obvious as $f$ is the ample model$/Z'$ of $K_Y+B_Y$.

(3) Since $K_{X'}+B'\sim_{\Rr,Z'}0$, the restriction of $K_{Y}+B_Y$ to a general fiber of the induced morphism $Y\rightarrow Z'$ is $\Rr$-linearly equivalent to $0$. Thus $Z$ is birational to $Z'$, so $\dim Y>\dim Z$. Since $q': \bar X'\rightarrow X'$ and the induced birational map $\bar X'\dashrightarrow Y$ are both isomorphisms over the generic point of $Z'$, $S_Y$ is horizontal$/Z'$, hence $S_Y$ is horizontal$/Z$.  This implies (3).

(4)(5) Since the induced birational map $\bar X\dashrightarrow Y$ is a sequence of steps of a $(K_{\bar X}+\bar B)$-MMP, we have that $Y$ is $\mathbb Q$-factorial and
$$\kappa_{\sigma}(K_{Y}+B_Y)=\kappa_{\sigma}(K_{\bar X}+\bar B)=\kappa_{\sigma}(K_X+B),$$
and
$$\kappa_{\iota}(K_{Y}+B_Y)=\kappa_{\iota}(K_{\bar X}+\bar B)=\kappa_{\iota}(K_X+B).$$
This implies (4). Moreoer, if $(X,B)$ is klt, then by our construction, $(\bar X,\bar B)$ is klt, so $(Y,B_Y)$ is $\mathbb Q$-factorial klt.
\end{proof}

\begin{thm}\label{thm: abundance lc case assuming non-vanishing}
Assume Conjecture \ref{conj: non-vanishing smooth} in dimension $d-1$ and for varieties of numerical dimension $1$ in dimension $d$. Then for any projective lc pair $(X,B)$ of dimension $d$, such that $0\leq\kappa_{\sigma}(K_X+B)\leq 1$, $(X,B)$ has a good minimal model. In particular, we have $\kappa_{\sigma}(K_X+B)=\kappa_{\iota}(K_X+B)$. 
\end{thm}
\begin{proof}
 By Theorem \ref{thm: 1.5 low imply 1.7}, Conjecture \ref{conj: stof with scaling weak} holds in dimension $d$. By Theorem \ref{thm: projective case}, we only need to show that $\kappa_{\iota}(K_X+B)\geq 0$. We may assume that Theorem \ref{thm: abundance lc case assuming non-vanishing} holds in dimension $\leq d-1$.

 Let $h: Y\rightarrow X$ be a log resolution of $(X,B)$ and write $B_Y:=h^{-1}_*B+\Exc(h)$. Let
$$t:=\inf\{s\geq 0\mid K_{Y}+sB_Y\text{ is pseudo-effective}\}.$$
Then we have 
$$1=\kappa_{\sigma}(K_X+B)=\kappa_{\sigma}(K_{Y}+B_Y)\geq\kappa_{\sigma}(K_{Y}+tB_Y)\geq 0$$
and
$$\kappa_{\iota}(K_X+B)=\kappa_{\iota}(K_{Y}+B_Y)\geq\kappa_{\iota}(K_{Y}+tB_Y).$$
Thus, possibly replacing $(X,B)$ with $(Y,tB_Y)$, we may assume that $(X,B)$ is log smooth, $h$ is the identity morphism, and either $B=0$ or $t=1$. By Theorem \ref{thm: num 0 abundance}, we may assume that $\kappa_{\sigma}(K_X+B)=1$. By our assumption we may assume that $B\not=0$. Let $S:=\lfloor B\rfloor$. There are two cases.

\medskip

\noindent\textbf{Case 1.} $K_X+B-\epsilon S$ is pseudo-effective for some $\epsilon>0$, then since 
$$1=\kappa_{\sigma}(K_X+B)\geq\kappa_{\sigma}(K_{X}+B-\epsilon S)\geq 0$$
and
$$\kappa_{\iota}(K_X+B)\geq\kappa_{\iota}(K_{Y}+B-\epsilon S),$$
possibly replacing $(X,B)$ with $(X,B-\epsilon S)$ for some $0<\epsilon\ll 1$, we may assume that $(X,B)$ is klt. Since $t=1$, there exists an irreducible component $T$ of $B$ such that $K_X+B-\delta T$ is not pseudo-effective for any $0<\delta\ll 1$. By Lemma \ref{lem: non pe induces contraction}, there exists a $\mathbb Q$-factorial projective klt pair $(Y,B_Y)$ of dimension $d$ and a contraction $f: Y\rightarrow Z$, such that 
\begin{itemize}
    \item $K_Y+B_Y\sim_{\mathbb R,Z}0$,
    \item $d=\dim Y>\dim Z$, and
    \item $\kappa_{\sigma}(K_Y+B_Y)=\kappa_{\sigma}(K_X+B)=1$ and $\kappa_{\iota}(K_Y+B_Y)=\kappa_{\iota}(K_X+B)$.
\end{itemize}
By Lemma \ref{lem: ambro05 r coefficient} there exists a klt pair $(Z,B_Z)$ such that
$$K_{Y}+B_Y\sim_{\mathbb R}f^*(K_Z+B_Z).$$
By induction hypothesis, 
$$\kappa_{\iota}(K_X+B)=\kappa_{\iota}(K_Y+B_Y)=\kappa_{\iota}(K_Z+B_Z)=\kappa_{\sigma}(K_Z+B_Z)=\kappa_{\sigma}(K_Y+B_Y)=1$$
and we are done.

\medskip

\noindent\textbf{Case 2.} $K_X+B-\epsilon S$ is not pseudo-effective for any $\epsilon>0$. Then there exists an irreducible component $T$ of $S$ such that $K_X+B-\epsilon T$ is not pseudo-effective for any $\epsilon>0$. By Lemma \ref{lem: non pe induces contraction}, there exists an lc pair $(Y,B_Y)$ of dimension $d$, a prime divisor $\tilde T_Y$ on $Y$ with normalization $T_Y$, and a contraction $f: Y\rightarrow Z$, such that
\begin{itemize}
    \item $\tilde T_Y$ is an irreducible component of $B_Y$ and $\mult_{\tilde T_Y}B_Y=1$.
    \item $K_Y+B_Y\sim_{\mathbb R,Z}0$,
    \item $d=\dim Y>\dim Z$ and $\tilde T_Y$ is horizontal$/Z$.
    \item $\kappa_{\sigma}(K_Y+B_Y)=\kappa_{\sigma}(K_X+B)=1$ and $\kappa_{\iota}(K_Y+B_Y)=\kappa_{\iota}(K_X+B)$.
\end{itemize}
Let
$$K_{T_Y}+B_{T_Y}:=(K_Y+B_Y)|_{T_Y}$$
and $g: T_Y\rightarrow Z$ the induced projective surjective morphism. Then $(T_Y,B_{T_Y})$ is a projective lc pair of dimension $d-1$ and there exists an $\mathbb R$-divisor $L$ on $Z$ such that
$$K_Y+B_Y=f^*L\text{ and }K_{T_Y}+B_{T_Y}=g^*L.$$
Thus
$$\kappa_{\sigma}(K_{T_Y}+B_{T_Y})=\kappa_{\sigma}(K_Y+B_Y)=1,$$
so by induction hypothesis,
$$\kappa_{\iota}(K_X+B)=\kappa_{\iota}(K_{Y}+B_{Y})=\kappa_{\iota}(K_{T_Y}+B_{T_Y})=1.$$
and we are done.
\end{proof}
\begin{proof}[Proof of Theorem \ref{thm: main nv}]
    This is a special case of Theorem \ref{thm: abundance lc case assuming non-vanishing}.
\end{proof}

\section{Applications and discussions}

\subsection{Further results on abundance}
By combining existing results on the non-vanishing conjecture with our main theorems, we obtain more new results on the abundance conjecture. For simplicity, we restrict our discussion to unconditional results (i.e., those not dependent on any major conjectures) in dimensions four and five.

\begin{cor}
    Let $X$ be a minimal projective terminal variety of dimension $\leq 5$ such that $\nu(X)=1$. Assume that there exists a positive integer $q$ such that
    $$h^0(X,\Omega_X^{[q]}\otimes\mathcal{O}_X(mK_X))>0$$
    for infinitely many $m$ such that $mK_X$ is Cartier. Then $K_X$ is semi-ample.
\end{cor}
\begin{proof}
    It follows from Theorem \ref{thm: lc abundance dim 5 nu 1} and \cite[Theorem A]{LP18}.
\end{proof}

\begin{cor}
Let $X$ be a minimal projective klt variety of dimension $\leq 5$ such that $\nu(X)=1$, and there exists $s\in H^0(X,\Omega_X^{[q]})$ which vanishes along some divisor. Then $K_X$ is semi-ample.
\end{cor}
\begin{proof}
        It follows from Theorem \ref{thm: lc abundance dim 5 nu 1} and \cite[Theorem 6.11]{LP18}.
\end{proof}

\begin{cor}
    Let $(X,B)$ be a projective lc pair in dimension $\leq 5$ such that $\kappa_{\sigma}(K_X+B)=1$ and $(X,\Delta)$ is klt Calabi-Yau for some $\Delta$. Then $(X,B)$ has a good minimal model. In particular, $\kappa_{\iota}(K_X+B)=1$.
\end{cor}
\begin{proof}
    It follows from Theorem \ref{thm: lc abundance dim 5 nu 1} and \cite[Theorem 1.3]{Has19b}.
\end{proof}

\begin{cor}
    Let $(X,B)$ be a projective uniruled lc pair in dimension $4$ such that $\kappa_{\sigma}(K_X+B)=1$. Then $(X,B)$ has a good minimal model. In particular, $\kappa_{\iota}(K_X+B)=1$.
\end{cor}
\begin{proof}
    It follows from Theorem \ref{thm: lc abundance dim 5 nu 1} and \cite[Corollary 1.2]{LM21}.
\end{proof}

\begin{cor}
    Let $(X,B)$ be a projective rationally connected lc pair of dimension $\leq 5$ such that $\kappa_{\sigma}(K_X+B)=1$. Then $(X,B)$ has a good minimal model. In particular, $\kappa_{\iota}(K_X+B)=1$.
\end{cor}
\begin{proof}
       It follows from Theorem \ref{thm: lc abundance dim 5 nu 1} and \cite[Theorem A.1]{JLX22}. 
\end{proof}

\subsection{The relative case}

As we mentioned in the introduction, the abundance conjecture can be reduced to the case when $\kappa=0$ or $-\infty$ thanks to \cite{GL13,Lai11}. Since we have obtained new results on the abundance conjecture when $\nu=1$, using similar arguments, we obtain the following result:

\begin{thm}\label{thm: relative abundance nu=1}
    Let $(X,B)/U$ be a klt pair such that $\dim X-\dim U=d$. Let $\kappa:=\kappa_{\iota}(X/U,K_X+B)$ and $\nu:=\kappa_{\sigma}(X/U,K_X+B)$. Assume that $\nu-\kappa\leq 1$ and one of the followings hold:
    \begin{enumerate}
    \item $d-\kappa\leq 5$.
    \item Conjecture \ref{conj: stof with scaling weak} holds in dimension $d-\kappa$.
    \item Conjecture \ref{conj: non-vanishing smooth} holds in dimension $d-\kappa-1$.
    \end{enumerate}
    Then $(X,B)/U$ has a good minimal model.
\end{thm}
\begin{proof}
First we prove the theorem when $U$ is a closed point. By \cite[Lemmas 3.6, 3.8]{Bir10} and Theorem \ref{thm: 1.5 low imply 1.7}, we only need to prove the theorem under the assumption of (2). Let $f: X'\rightarrow Z$ be an invariant Iitaka fibration (cf. \cite[Definition 6.1]{CHL24}) of $K_X+B$. Possibly replacing $X'$ and $Z$ with higher resolutions, we may assume that $Z$ is smooth and the induced birational map $h: X'\dashrightarrow X$ is a log resolution of $(X,B)$. We let $B':=h^{-1}_*B+(1-\delta)\Exc(h)$ for some $0<\delta\ll 1$. Then
$$K_{X'}+B'=h^*(K_X+B)+E$$
for some $E\geq 0$ that is exceptional$/X$. Since $h$ is the invariant Iitaka fibration of $K_X+B$, we have
$$\dim Z=\kappa,$$
and there exists an ample ample$/U$ $\mathbb R$-divisor $A\geq 0$ on $Z$ and an exceptional$/X$ $\mathbb R$-divisor $E'\geq 0$ on $X'$, such that 
$$h^*(K_X+B)\sim_{\mathbb R}f^*A+E'.$$
Then for any real number $k$, we have
$$K_{X'}+B'+kf^*A\sim_{\mathbb R}(1+k)f^*A+E+E'.$$
By Lemmas \ref{lem: property of numerical and Iitaka dimension} and \ref{lem: compare kappa}, for any $k\geq 0$ we have 
\begin{align*}
    \kappa_{\sigma}(K_{X'}+B'+kf^*A)&=\kappa_{\sigma}((1+k)f^*A+E+E')=\kappa_{\sigma}(f^*A+E+E')\\
    &=\kappa_{\sigma}(h^*(K_X+B)+E)=\kappa_{\sigma}(K_{X'}+B')=\kappa_{\sigma}(K_X+B)=\nu.
\end{align*}
By \cite[(3.3)]{Fuj20}, we have
 \begin{align*}
\nu &=    \kappa_{\sigma}(K_{X'}+B'+kf^*A)\geq \kappa_{\sigma}(X'/Z,K_{X'}+B')+\kappa(A)\\
    &=\kappa_{\sigma}(X'/Z,K_{X'}+B')+\dim Z=\kappa_{\sigma}(X'/Z,K_{X'}+B')+\kappa.
\end{align*}
Therefore, $\kappa_{\sigma}(X'/Z,K_{X'}+B')\in\{0,1\}$.  

Let $F$ be a general fiber of $f$ and let $K_F+B_F:=(K_{X'}+B')|_F$. Then $\kappa_{\iota}(K_F+B_F)=0$ and $\kappa_{\sigma}(K_F+B_F)\in\{0,1\}$. We have
$$\dim F=\dim X-\dim Z=d-\kappa.$$
By Theorem \ref{thm: projective case}, $\kappa_{\sigma}(K_F+B_F)=0$, so $(F,B_F)$ has a good minimal model. By \cite[Theorem 1.2]{Has19a}, $(X',B')/Z$ has a good minimal model $(X'',B'')/Z$, and we have $K_{X''}+B''\sim_{\mathbb R,Z}0$. By Lemma \ref{lem: ambro05 r coefficient}, there exists a klt pair $(Z,B_Z)$ such that
$$K_{X''}+B''\sim_{\mathbb R}f''^*(K_{Z}+B_Z)$$
where $f'': X''\rightarrow Z$ is the induced contraction. Then
$$\kappa_{\iota}(K_{Z}+B_Z)=\kappa_{\iota}(K_{X''}+B'')=\kappa_{\iota}(K_{X'}+B')=\kappa=\dim Z,$$
so $K_Z+B_Z$ is big, hence 
$$\kappa_{\sigma}(K_{X''}+B'')=\kappa_{\sigma}(K_{Z}+B_Z)=\dim Z.$$ 
By \cite[Theorem 1.2]{BCHM10}, $(X'',B'')$ has a good minimal model, hence $(X',B')$ has a good minimal model. By Lemma \ref{lem: equivalence of log minimal model}, $(X,B)$ has a good minimal model.

Now we prove the theorem when $U$ is not necessarily a closed point. Possibly passing to the Stein factorization, we may assume that $X\rightarrow U$ is a contraction. Let $G$ be a general fiber of $X\rightarrow U$ and let $K_G+B_G:=(K_X+B)|_G$. By the projective case we have just proved, $(G,B_G)$ has a good minimal model. By \cite[Theorem 1.2]{Has19a}, $(X,B)/U$ has a good minimal model.
\end{proof}

\begin{rem}\label{rem: why lc not ok}
At the moment, in Theorem \ref{thm: relative abundance nu=1}, we have to assume that $(X,B)$ is klt. This is because two technical difficulties. The first difficulty is that results on the existence of good minimal models usually require that ``all lc centers dominate the base", hence existence of good minimal models along a general fiber usually could not extend to a relative good minimal model. The second difficulty is that we do not have the existence of good minimal models for lc general type varieties and even the $\kappa=\nu=d$ case could not be proven. We still expect more partial results to hold if some log abundant properties are assumed.
\end{rem}

\subsection{Further discussions}

We thank Professor Vladimir Lazić for suggesting the following remark.

\begin{rem}\label{rem: kaw11's idea}
A note by Professor Kawamata \cite{Kaw11} claimed that the abundance conjecture in dimension $d$ follows from the following three conditions:
\begin{enumerate}
\item The abundance conjecture in dimensions $\leq d-1$.
\item The (log) non-vanishing conjecture in dimension $d$.
\item The termination of flips in dimension $d$.
\end{enumerate}
However, \cite{Kaw11} was withdrawn after four days ``due to a crucial error”. The key issue lies in the fact that, under the setting of Theorem \ref{thm: special abundance}, the assumption that the numerical dimension is one in condition (3) is essential in order to deduce $K_S + B_S \equiv 0$, which is in turn a critical step in the proof of Theorem \ref{thm: special abundance}. In other words, Kawamata’s approach remains valid in the case $\nu = 1$, which follows essentially the same ideas as in \cite{Kaw92} combining with the use of the slc property of Du Bois singularities \cite{KK10,Kol13}. Thus, despite the withdrawal of \cite{Kaw11}, the conceptual contributions it introduced, some of which are adopted in this paper, should not be overlooked.

On the other hand, even when restricting to the numerical dimension one case, the approach in \cite{Kaw11} does not imply any of our main theorems, as the assumptions it requires are too strong. A key reason we are able to remove or weaken these assumptions is the significant progress in the minimal model program in recent years. Many of these advances, e.g. \cite{Has18,HH20,LT22,MZ23}, are highly nontrivial and rely essentially on the ACC for log canonical thresholds, as explained in \cite{Bir07}. We note that the ACC for lc thresholds was only proven in 2012 \cite[Theorem 1.1]{HMX14}, well after the appearance of \cite{Kaw11}. Our paper is made possible thanks to these recent developments.
\end{rem}

\end{document}